\documentclass[letterpaper, preprint, paper,11pt]{AAS}	

\usepackage{bm}
\usepackage{amsmath}
\usepackage{amssymb}
\usepackage{subfig}
\usepackage[colorlinks=true, pdfstartview=FitV, linkcolor=black, citecolor= black, urlcolor= black]{hyperref}
\usepackage{overcite}
\usepackage{footnpag}			      	
\usepackage{comment}

\newcommand{\hide}[1]{}  

\newcommand{\STMa}[2]{\bm{\Phi}^{\textbf{#1}}_{\textbf{#2}}}

\newcommand{\STMb}[2]{\bm{\Phi}^{\textbf{#1}}_{\bm{#2}}}

\newcommand{\STMc}[2]{\bm{\Phi}^{\bm{#1}}_{\textbf{#2}}}

\newcommand{\STMd}[2]{\bm{\Phi}^{\bm{#1}}_{\bm{#2}}}

\PaperNumber{24-206}

\begin{document}


\title{Generation of Energy-Optimal Low-Thrust forced periodic trajectories in the CR3BP}


\author{Colby C. Merrill\thanks{PhD Candidate, Sibley School of Mechanical and Aerospace Engineering, Cornell University, Ithaca, NY 14853, USA}, 
Jackson Kulik\thanks{Assistant Professor, Department of Mechanical and Aerospace Engineering, Utah State University, Logan, UT 84322, USA}, 
Dmitry Savransky\thanks{Associate Professor, Sibley School of Mechanical and Aerospace Engineering, Cornell University, Ithaca, NY 14853, USA}}

\maketitle{}


\begin{abstract} 
In this work, we investigate trajectories that require thrust to maintain periodic structure in the circular restricted three-body problem (CR3BP). We produce bounds in position and velocity space for the energy-constrained reachable set of initial conditions. Our trajectories are energy-optimal and analyzed via linear analysis. We provide validation for our technique and analyze the cost of deviating in various directions to the reference. For our given reference, we find that it is relatively expensive to decrease perilune distance for orbits in the Earth-Moon system.
\end{abstract}


\section{Introduction}

Uncontrolled, periodic orbits in the circular restricted three-body problem (CR3BP) have been well-studied and are often exploited in spacecraft operations. There exist families of these periodic orbits centered around the five equilibrium Lagrange points in a three-body system. Of particular interest are periodic orbits around the L1 and L2 points, as these orbits have historically been used by spacecraft in the Earth-Sun system and are planned for use in the Earth-Moon system.

A constant perturbation applied to the CR3BP (e.g., a constant acceleration provided by low thrust propulsion), will shift the equilibrium points and allow areas beyond the natural space to be exploited \cite{Cox2019, Cox2020}. Because the equilibrium points have been shifted, the periodic trajectories around the equilibria will be shifted as well \cite{Morimoto2006, Tsuruta2024}. Many previous studies have primarily focused on constant, low thrust trajectories \cite{Cox2019, Cox2020, Morimoto2006, DeLeo2022} or optimal control in proximity of the shifted equilibria \cite{Tsuruta2024}. Other related work exists that focuses on forced circumnavigation or controlled loitering trajectories relative to a spacecraft on some reference orbit in the CR3BP \cite{Sandel2024}. In this work, we also investigate optimally controlled trajectories relative to a reference orbit. However, the emphasis here is on developing and characterizing a catalog of forced periodic trajectories (trajectories that begin and end at the same state) under control rather than on the relative motion of two real satellites. Additionally, the lens through which we examine optimally controlled trajectories in the vicinity of a reference orbit leverages the analytical aspects of techniques from relative reachable set theory to provide a lower fidelity but more holistic view of optimally controlled trajectories near a reference orbit. We refer to these orbits as forced periodic trajectories. While not periodic in the natural dynamics, these trajectories expand available options for satellite orbits in the CR3BP beyond naturally periodic or quasi-periodic structures.


\section{Methods}
\label{s:methods}

As an initial investigation into this problem, we make use of the linearized model of energy-optimal control and approximate the energy cost of optimal forced periodic trajectories near a naturally periodic orbit about the L2 point in the Earth-Moon system. This periodic orbit will serve as a reference trajectory for our spacecraft and we will work to understand how small perturbations in the initial and final state from this reference trajectory affect the cost of control. 

\subsection{The Energy-Limited Reachable Set}
Here, we are interested in finding the reachable set of states for our system in the energy-constrained problem. This formalism is based on the work in Kulik et al., 2024 \cite{Kulik2024}, where it is expanded upon further. We specialize results from this work to the context of forced periodic trajectories. We assume a dynamical system of the form
\begin{equation}
    \frac{\mathrm{d}\mathbf{x}}{\mathrm{d}t}=\mathbf{F}(\mathbf{x})+\begin{bmatrix}\mathbf{0}_{3\times1}\\
        \mathbf{u}
    \end{bmatrix}
\end{equation}
where the state vector $\mathbf{x}\in\mathbb{R}^6$ is defined by stacking the position and velocity vectors $\mathbf{x}=[\mathbf{r}^T, \mathbf{v}^T]^T$ and $\mathbf{u}$ is the control acceleration vector. In our system, $\mathbf{F}(\mathbf{x})$ gives the natural dynamics for the third body in the ``canonical rotating frame" of the CR3BP
\begin{equation}
    \mathbf{F}(\mathbf{x}) =
    \begin{bmatrix}
        v_{x}\\ v_{y}\\ v_{z}\\ 
        2v_{y} + x - (1-\mu^*)\frac{x+\mu^*}{R_1^3} - \mu^*\frac{x-1+\mu^*}{R_2^3}\\
        -2v_{x} + y - (1-\mu^*)\frac{y}{R_1^3} - \mu^*\frac{y}{R_2^3}\\
        -(1-\mu^*)\frac{z}{R_1^3} - \mu^*\frac{z}{R_2^3}
    \end{bmatrix}
\end{equation}
where $x$, $y$, and $z$ are the components of the spacecraft's position vector, $\mathbf{r}$, and $v_x$, $v_y$, and $v_z$ are the components of the spacecraft's velocity vector, $\mathbf{v}$. 
The distances of the third body with respect to the primary and secondary are defined as $R_1$ and $R_2$, respectively, and are evaluated as
\begin{align}
    R_1 &= \sqrt{(x+\mu^*)^2 + y^2 + z^2}\\
    R_2 &= \sqrt{(x-1+\mu^*)^2 + y^2 + z^2}
\end{align}
where $\mu^* = m_2/(m_1+m_2)$ is the mass parameter of the system. Defining a quadratic cost function of the form
\begin{equation}
    J = \frac{1}{2}\int_{t_0}^{t_f} ||\mathbf{u}||^2 dt
\end{equation}
the optimal control from one state to another is given by solving a two-point boundary value problem associated with a system of ordinary differential equations. These equations have twice as many dimensions as the state of the original system and are given by
\begin{align}
    \frac{\mathrm{d}\mathbf{x}}{\mathrm{d}t}&=\mathbf{F}(\mathbf{x})+\begin{bmatrix}
        \mathbf{0}_{3\times1}\\\mathbf{u}\end{bmatrix}\\
    \frac{\mathrm{d}\bm{\lambda}}{\mathrm{d}t}&=-\left(\frac{\partial \mathbf{F(x)}}{\partial \mathbf{x}}\right)^T \bm{\lambda}    \\
    \mathbf{u} &= -\bm{\lambda_v}
\end{align}
where $\bm{\lambda_v}$ is the velocity costate vector given by the last three elements of the costate vector \cite{bryson2018applied}. $J$ can then be written in terms of the velocity costate vector as
\begin{equation}
    J = \frac{1}{2}\int_{t_0}^{t_f} \bm{\lambda_v}^T \bm{\lambda_v} dt
\end{equation}
The augmented state is given by stacking the states and costates
\begin{equation}
\mathbf{y} = \begin{bmatrix}
        \mathbf{x} \\ 
        \bm{\lambda}
    \end{bmatrix} =
    \begin{bmatrix}
        \mathbf{r} \\
        \mathbf{v} \\
        \bm{\lambda}_r \\
        \bm{\lambda_v}
    \end{bmatrix}
\end{equation}
where $\mathbf{r}$ is the position vector from the origin to the spacecraft, $\mathbf{v}$ is the velocity of the spacecraft, and $\bm{\lambda}$ are the six costates. The state transition matrix (STM) associated with the augmented state vector and its dynamics yields a linear approximation of perturbations to the final augmented state $\delta\mathbf{y}(t_f)$ at some final time, $t_f$, as a function of deviations in the initial augmented state:
\begin{equation}
\delta\mathbf{y}(t_f) = \begin{bmatrix}
        \delta\mathbf{r}(t_f) \\ \delta\mathbf{v}(t_f) \\ \delta\bm{\lambda}_r(t_f) \\ \delta\bm{\lambda_v}(t_f)
    \end{bmatrix} \approx \STMd{}{}(t_f,t_0) \delta\mathbf{y}_0 = \begin{bmatrix}
      \STMa{x}{x} & \STMb{x}{\lambda}  \\ \STMc{\lambda}{x} & \STMd{\lambda}{\lambda}
    \end{bmatrix} \delta\mathbf{y}_0 = \begin{bmatrix}
      \STMa{r}{r} & \STMa{r}{v} & \STMb{r}{\lambda_r} & \STMb{r}{\lambda_v} \\ \STMa{v}{r} & \STMa{v}{v} & \STMb{v}{\lambda_r} & \STMb{v}{\lambda_v} \\ \STMc{\lambda_r}{r} & \STMc{\lambda_r}{v} & \STMd{\lambda_r}{\lambda_r} & \STMd{\lambda_r}{\lambda_v} \\ \STMc{\lambda_v}{r} & \STMc{\lambda_v}{v} & \STMd{\lambda_v}{\lambda_r} & \STMd{\lambda_v}{\lambda_v}
    \end{bmatrix} \delta\mathbf{y}_0
\end{equation}
where $\delta\mathbf{y}_0 = \delta\mathbf{y}(t_0)$ is the perturbed initial state and $\bm{\Phi}(t_f,t_0)$ is a time-varying STM associated with the augmented state, reference trajectory, and initial and final times. We adopt the notation
\begin{equation}
    \STMa{b}{a}(t,t_0) = \frac{\partial \mathbf{b}(t)}{\partial \mathbf{a}(t_0)}
\end{equation}
If the time dependence of an STM is omitted in this paper, it indicates that the STM corresponds to a full period (i.e., $\STMd{}{}(t_f,t_0) = \STMd{}{}$). $\delta\bm{\lambda_v}(t)$ at all times may be evaluated as 
\begin{equation}
\delta\bm{\lambda_v}(t) = \STMc{\lambda_v}{y}(t,t_0) \delta\mathbf{y}_0 = 
    \begin{bmatrix}
      \STMc{\lambda_v}{r}(t,t_0) & \STMc{\lambda_v}{v}(t,t_0) & \STMd{\lambda_v}{\lambda_r}(t,t_0) & \STMd{\lambda_v}{\lambda_v}(t,t_0)
    \end{bmatrix} \delta\mathbf{y}_0
\end{equation}
Substituting in to the cost function, we then have 
\begin{equation}
    J = \frac{1}{2}\delta\mathbf{y}_0^T \left(\int_{t_0}^{t_f} \left(\STMc{\lambda_v}{y}(t,t_0)\right)^T \left(\STMc{\lambda_v}{y}(t,t_0)\right) dt\right) \delta\mathbf{y}_0
    \label{eq:dy_cost_integral}
\end{equation}
With the linearized analysis, we explicitly solve for the initial costates that satisfy the linearized boundary value problem constraints
\begin{equation}
    \delta\bm{\lambda}_0 = \begin{bmatrix}
    -\left(\STMb{x}{\lambda}\right)^{-1} \STMa{x}{x} & \left(\STMb{x}{\lambda}\right)^{-1}  \end{bmatrix} \begin{bmatrix}
    \delta\mathbf{x}_0 \\ \delta\mathbf{x}_f
    \end{bmatrix}
\end{equation}
which then can be used to find the full, initial, augmented state in terms of the boundary conditions
\begin{equation}
    \delta\mathbf{y}_0 = \begin{bmatrix}
    \mathbf{I}_6 & \mathbf{0}_6 \\ -\left(\STMb{x}{\lambda}\right)^{-1} \STMa{x}{x} & \left(\STMb{x}{\lambda}\right)^{-1} \end{bmatrix} \begin{bmatrix}
    \delta\mathbf{x}_0 \\ \delta\mathbf{x}_f
    \end{bmatrix}
\end{equation}
We define the matrix
\begin{equation}
    \mathbf{E} =  \begin{bmatrix} 
    \mathbf{I}_6 & \mathbf{0}_6 \\ -\left(\STMb{x}{\lambda}\right)^{-1} \STMa{x}{x} & \left(\STMb{x}{\lambda}\right)^{-1} 
    \end{bmatrix}^T \int_{t_0}^{t_f} \left(\STMc{\lambda_v}{y}(t,t_0)\right)^T \left(\STMc{\lambda_v}{y}(t,t_0)\right) dt 
    \begin{bmatrix} \mathbf{I}_6 & \mathbf{0}_6 \\ -\left(\STMb{x}{\lambda}\right)^{-1} \STMa{x}{x} & \left(\STMb{x}{\lambda}\right)^{-1} 
    \end{bmatrix}
\end{equation}
which can be used to determine the energy-constrained reachable set for our system \cite{lee2018reachable, sun2020analysis}. Substituting in to the cost function, we now have 
\begin{equation}
    J = \frac{1}{2}\begin{bmatrix}
    \delta\mathbf{x}_0 \\ \delta\mathbf{x}_f
    \end{bmatrix}^T \mathbf{E}
    \begin{bmatrix}
    \delta\mathbf{x}_0 \\ \delta\mathbf{x}_f
    \end{bmatrix}
\end{equation}
Relative reachable set theory considers the set of $\delta\mathbf{x}_f$ that can be reached given some fixed starting condition $\delta\mathbf{x}_0$ or some set of potential starting points $\delta\mathbf{x}_0$. Much of this energy-limited reachable set theory is accomplished using ellipsoidal geometry associated with the quadratic form described by the matrix $\mathbf{E}$ and the eigenvalue decomposition. We will use these same techniques as well as extensions to generalized eigenvalues and eigenvectors to study the set of forced periodic trajectories that satisfy an energy constraint.

Forced periodic trajectories in the vicinity of a periodic reference orbit satisfy the condition $\delta\mathbf{x}_f=\delta\mathbf{x}_0$, so that the final state is equivalent to the initial state after one period of the reference periodic orbit. To study the set of forced periodic trajectories that require less than some energy limit, we may study the following matrix 
\begin{equation}
    \mathbf{E}^* = \begin{bmatrix}
        \mathbf{I}_6 & \mathbf{I}_6
    \end{bmatrix} \mathbf{E} \begin{bmatrix}
        \mathbf{I}_6 & \mathbf{I}_6
    \end{bmatrix}^T
\end{equation}
so that the linearized cost function for the forced periodic trajectory beginning and ending at $\delta\mathbf{x}_0$ can then be written as
\begin{equation}
    J = \frac{1}{2}\delta\mathbf{x}_0^T \mathbf{E}^* \delta\mathbf{x}_0
\end{equation}
Now assume that ($\gamma_i, \mathbf{w}_i$) is an eigenpair of the matrix $\mathbf{E}^*$ where $\gamma_i$ is some eigenvalue and $\mathbf{w}_i$ is its corresponding eigenvector. Since $\mathbf{E}^*$ is a symmetric positive semi-definite matrix, the set of possible relative states $\delta\mathbf{x}_0$ that cost less than some energy limit $J^*$ to begin and end at under linearized optimal control is given by the hyperellipsoid described by the set
\begin{equation}
    \left\{ \delta\mathbf{x}_0 \quad \mathrm{s.t.}\quad \frac{1}{2}\delta\mathbf{x}_0^T \mathbf{E}^* \delta\mathbf{x}_0 \leq J^* \right\}
\end{equation}
whose semi-axes are given by the vectors
\begin{equation}
    \mathbf{a}_i=\sqrt{\frac{2J^*}{\gamma_i}}\mathbf{w}_i
    \label{eq:semiaxes}
\end{equation}

Thus, we have defined the set of relative states which can be returned to in a period of the reference orbit with cost no greater than $J^*$. The ellipsoid described here is in 6-dimensional position and velocity space.


\section{Investigation and Analysis}

Our initial investigation involves selecting a reference halo orbit around L2. Because the reference orbit is naturally periodic, the initial phasing may be chosen at will. Here, we select the initial condition near apolune. The initial conditions for our chosen orbit are 
\begin{equation}
\begin{bmatrix}
    x\\ y\\ z\\ v_x\\ v_y\\ v_z\\ \bm{\lambda}
\end{bmatrix} = 
    \begin{bmatrix}
    1.06315768  \:\: \text{DU}\\  0.000326952322 \:\: \text{DU}\\ -0.200259761 \:\: \text{DU}\\  0.000361619362 \:\: \text{DU/TU}\\ -0.176727245 \:\: \text{DU/TU}\\ -0.000739327422 \:\: \text{DU/TU}\\ \mathbf{0}_{6\times1}
    \end{bmatrix} 
\end{equation}
with a period of 2.085034838884136 TU and a mass constant of 0.01215059. 

Propagation of our states and STMs is performed numerically in canonical units with an explicit 8th order Runge-Kutta integrator with absolute and relative error tolerances set to 1e-13. All computation is handled through precomputed STMs for the reference trajectory, as used in previous studies \cite{Kulik2023}.

\subsection{Validation}

In order to validate our method, we use a Newton-Raphson iteration scheme to find the true energy-optimal cost of our target orbits. We first determine the ``inherent cost'' to maintain the periodic reference orbits that we select. Although in theory these naturally periodic orbits should require no cost to maintain, there is typically some inherent numerical error in the initial conditions. Thus, the ``true'' cost includes error from the imperfect initial conditions. This is important for the validation analysis because the cost found by the Newton-Raphson scheme is inaccurate for cases where the cost would be lower or near the orbit's inherent cost. Therefore, we only consider the costs of orbits that are at least an order of magnitude greater than this inherent cost. Our reference trajectory's inherent cost is $\approx$ 3.5e-14 DU$^2$/TU$^3$. 

\begin{figure}[ht]
    \centering
    \subfloat{\includegraphics[width=0.5\textwidth]{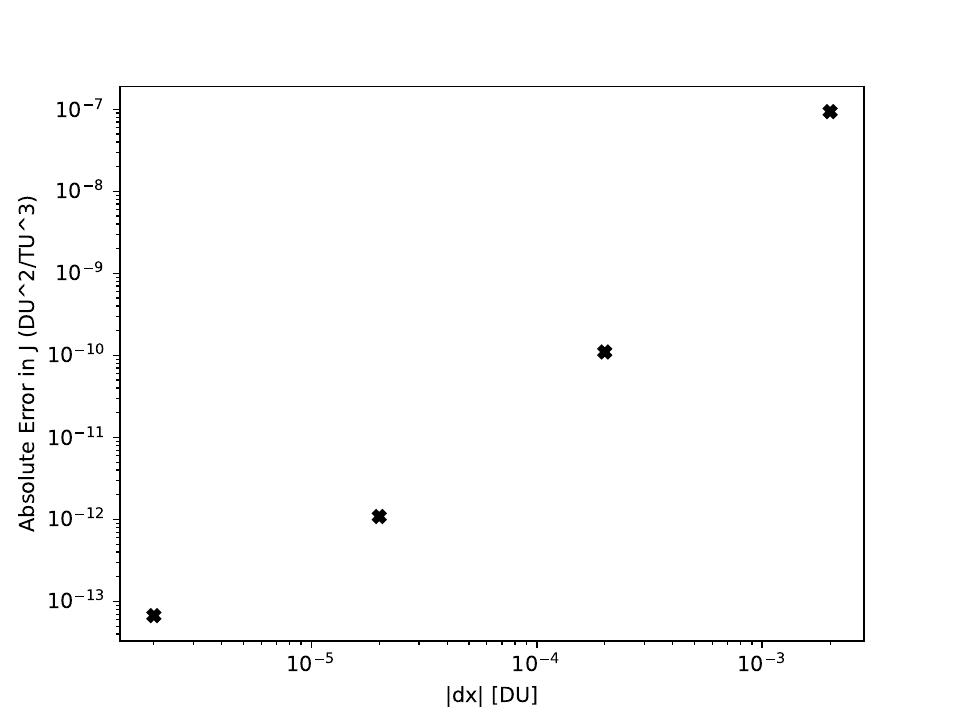}\label{fig:f1}}
  \subfloat{\includegraphics[width=0.5\textwidth]{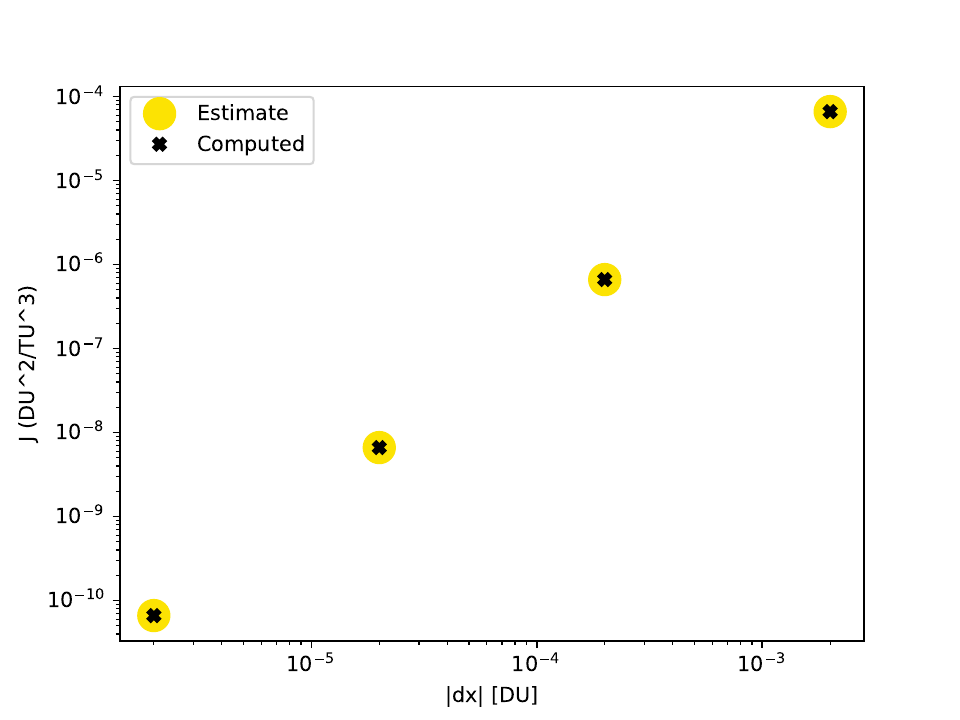}\label{fig:f2}}
  \caption{The left panel shows the difference between the estimated and computed costs while the right panel shows the estimated and computed costs. These deviations are made in the direction of the 4th eigenvector in Table \ref{tab:eigen}.}
    \label{fig:costs}
\end{figure}

In Figure \ref{fig:costs}, we see that the cost to maintain a periodic orbit follows a quadratic form and that the linear analysis provides a good approximation of the full nonlinear analysis. In the range of Figure \ref{fig:costs}, the error is at least 3 orders of magnitude below the true cost value (corresponding to $<$ 0.1\% error). The error grows as $\delta \mathbf{x}_0$ increases, as expected in this analysis. For deviations less than that shown in Figure \ref{fig:costs}, the linear analysis should be a good approximation of the true cost, although the inherent cost does begin to dominate the error.

\subsection{Addition of Thrust Constraint}

We now introduce a thrust constraint to analyze the extents of the energy-limited reachable set. For a given spacecraft with maximum acceleration from thrust $u_\text{max}$ and with the assumption that the spacecraft will thrust for the entire orbital period, we can find 
\begin{equation}
    J^* = \frac{1}{2} u_{\text{max}}^2 (t_f-t_0)
\end{equation}
where $t_f-t_0$ is the time required to complete a single period. As it is written, $J^*$ may be understood as the maximum specific power expended during a single orbit. With this constraint, we can solve for the extents of the semiaxes of the energy-limited reachable set. For a base design, let us assume a representative spacecraft with an initial mass of 1000 kg and thrust of 50 mN, translating to $u_\text{max}$ = 5E-5 m/s$^2 \: \approx$ 0.0184 DU/TU$^2$ and $J^*$ = 3.51E-4 DU$^2$/TU$^3$. This cost constraint is selected, in particular, because it falls in the range of costs with $<$ 1\% error and is a currently achievable thrust level for a propulsion system. Thus, it makes for a good example of this analysis.

\begin{table}[h]
    \centering
    \caption{The extents and corresponding directions of the semiaxes for the cost, $J^*$ = 3.51E-4 DU$^2$/TU$^3$. The directions are given by the eigenvectors of $\mathbf{E}^*$. The extents are listed from largest to smallest, corresponding to least expensive to most expensive directions to deviate from the reference in.}
    \begin{tabular}{c c}
    \hline
    Extent [DU] & Direction [DU, DU, DU, DU/TU, DU/TU, DU/TU]\\
    \hline
        15918.250 & [0.00075547 -0.36919333 -0.00154457 -0.40833012 -0.00219097  0.83483833] \\
 0.01984495 & [-0.4073362  -0.00359944 -0.29024187 -0.00429968  0.86591161 -0.00159068] \\
 0.00892213 & [ 0.00166912 -0.0949137   0.00474959 -0.87702446 -0.00323734 -0.47093912] \\
 0.00460856 & [ 0.58703528 -0.00312857  0.64311054  0.00223334  0.4917122   0.00165787] \\
 0.00244912 & [ 0.69940335  0.02398051 -0.70835806 -0.00809316  0.09164505  0.00494355] \\
 0.00051309 & [-0.01727422  0.92416983  0.01929803 -0.25299333  0.00145136  0.28501156] \\
    \end{tabular}
    \label{tab:eigen}
\end{table}

Using the thrust constraint paired with the linear analysis, we can study the reachable set. The reachable set takes the shape of a 6-dimensional hyper-ellipsoid with extents and semi-axes directions specified by the thrust constraint and eigenstructure of $\mathbf{E}^*$, as listed in Table \ref{tab:eigen}. The directions of the semi-axes are given by the eigenvectors of $\mathbf{E}^*$ and the extents of the semi-axes can be found by taking the magnitude of Equation \ref{eq:semiaxes}. The eigenvector associated with the smallest eigenvalue corresponds to the least expensive direction to deviate from the reference trajectory, which presents itself as the largest semi-axis extent. In this case, the largest extent should be infinite, as it corresponds to the zero eigenvalue. However, there is numerical error and so the first eigenvalue has a finite length here. For this reason, we do not present results that correspond to deviations in the direction of the first eigenvector.

\begin{figure}[ht]
    \centering
    \subfloat{\includegraphics[width=0.5\textwidth]{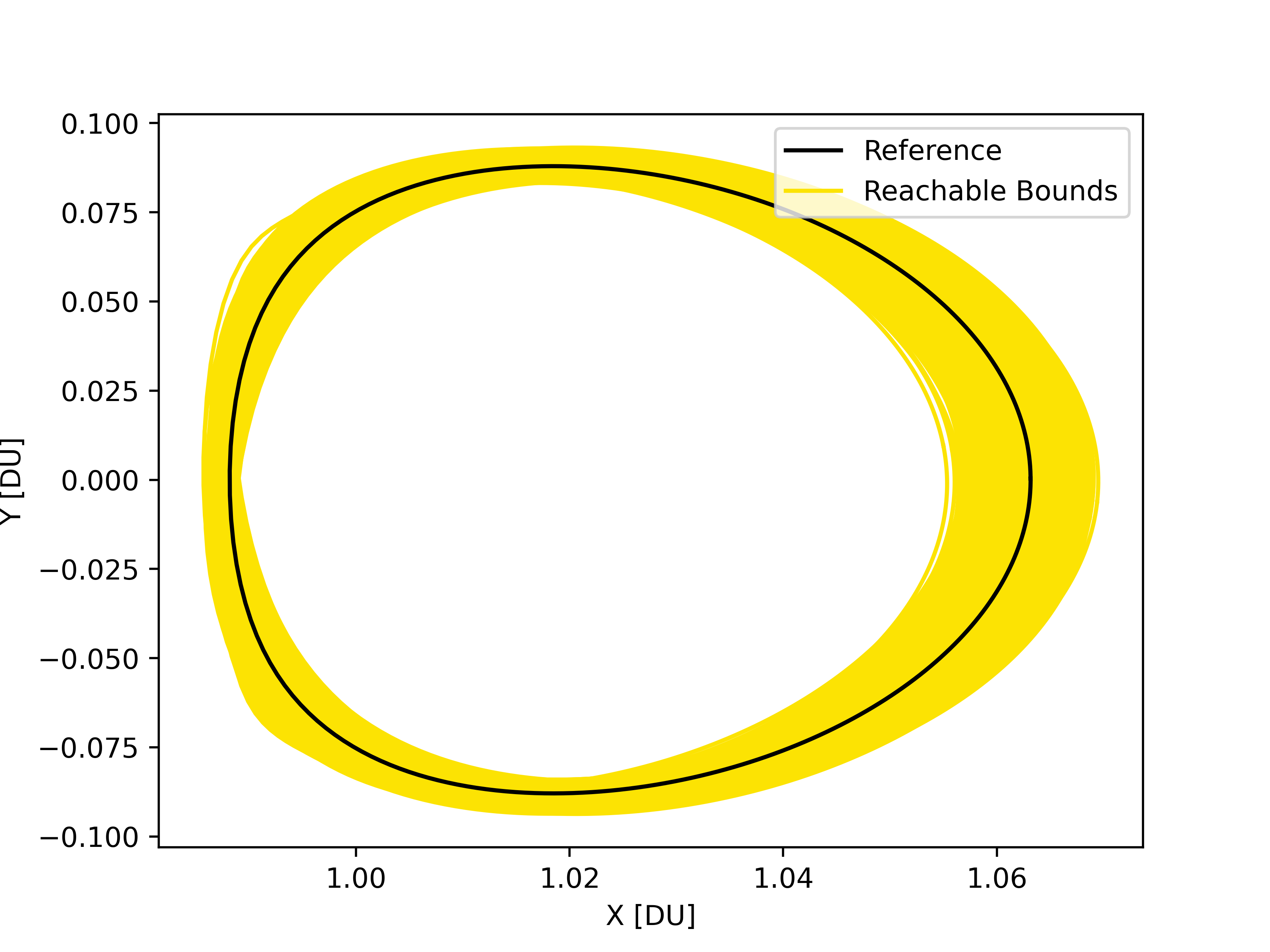}}
  \subfloat{\includegraphics[width=0.5\textwidth]{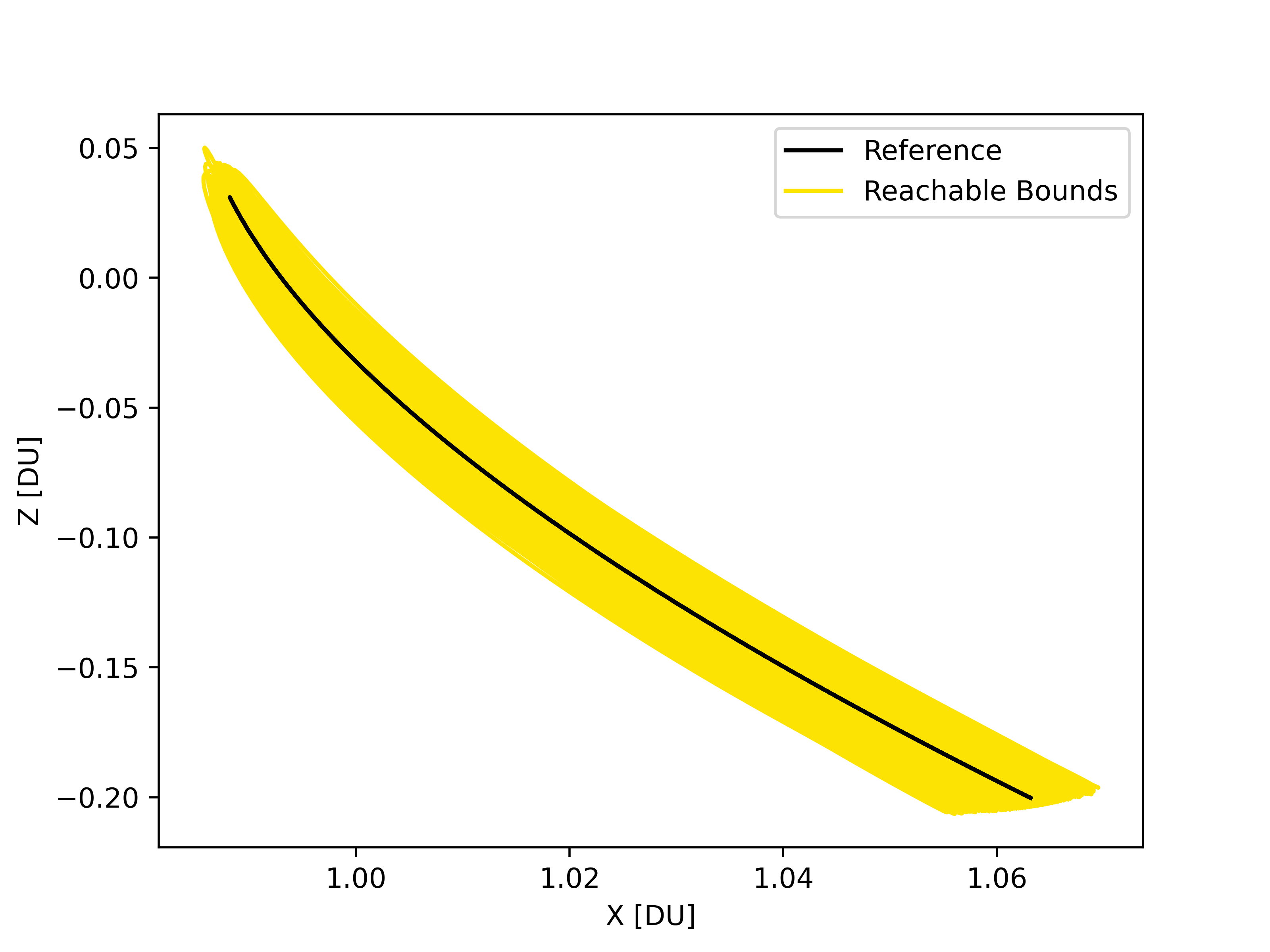}}
  \\
  \subfloat{\includegraphics[width=0.5\textwidth]{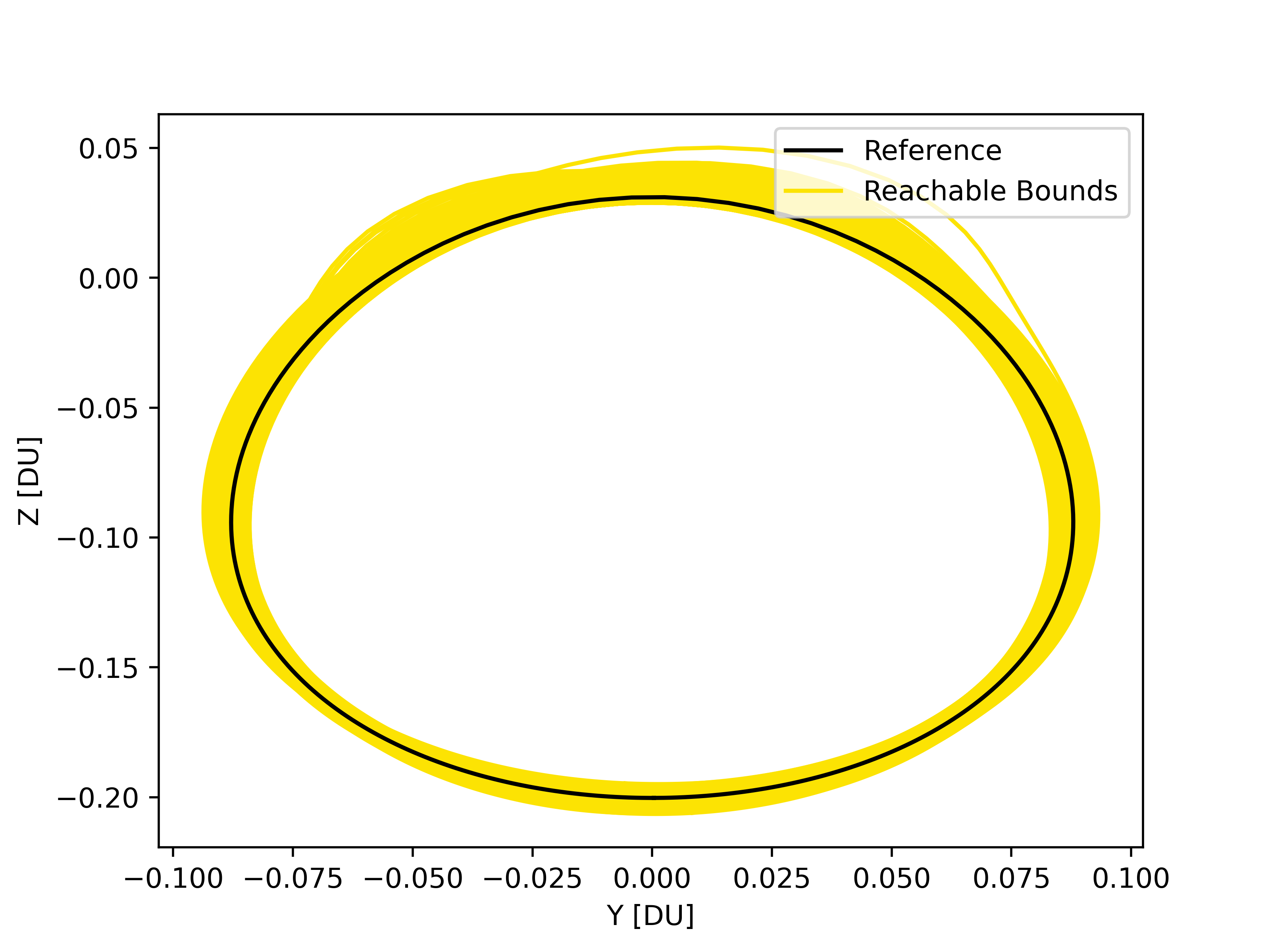}}
  \subfloat{\includegraphics[width=0.5\textwidth,clip=true, trim=50pt 0pt 20pt 50pt]{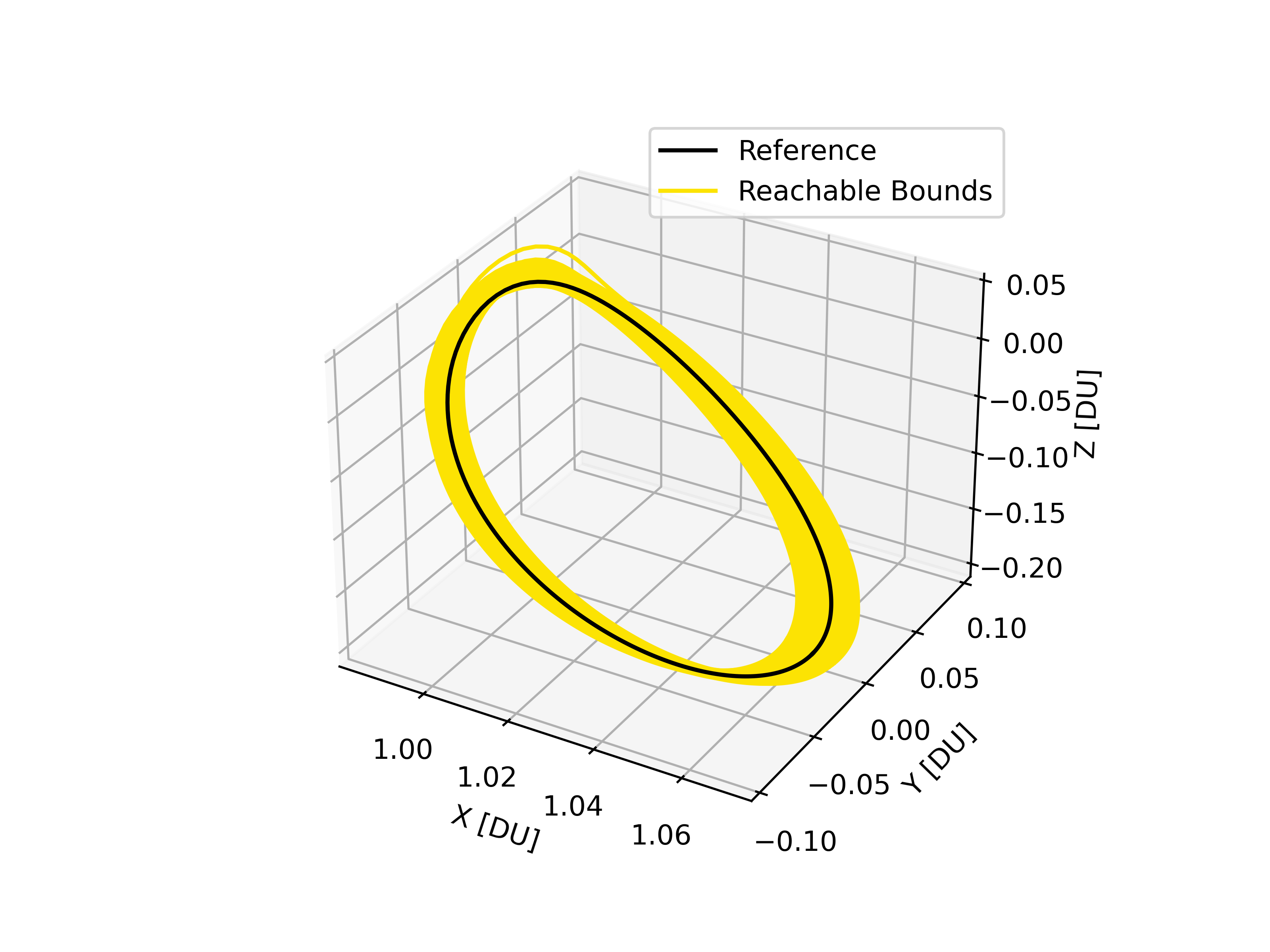}}
  \caption{The reachable sets in position space relative to the reference trajectory via sampling. The initial conditions are sampled 100,000 times and propagated for a period to produce these orbits. The direction of the initial deviation is randomized but the cost for each is $J^*$ = 3.51E-4 DU$^2$/TU$^3$.}
    \label{fig:Trajectories}
\end{figure}

\begin{figure}[htbp]
    \centering
    \subfloat{\includegraphics[width=0.5\textwidth]{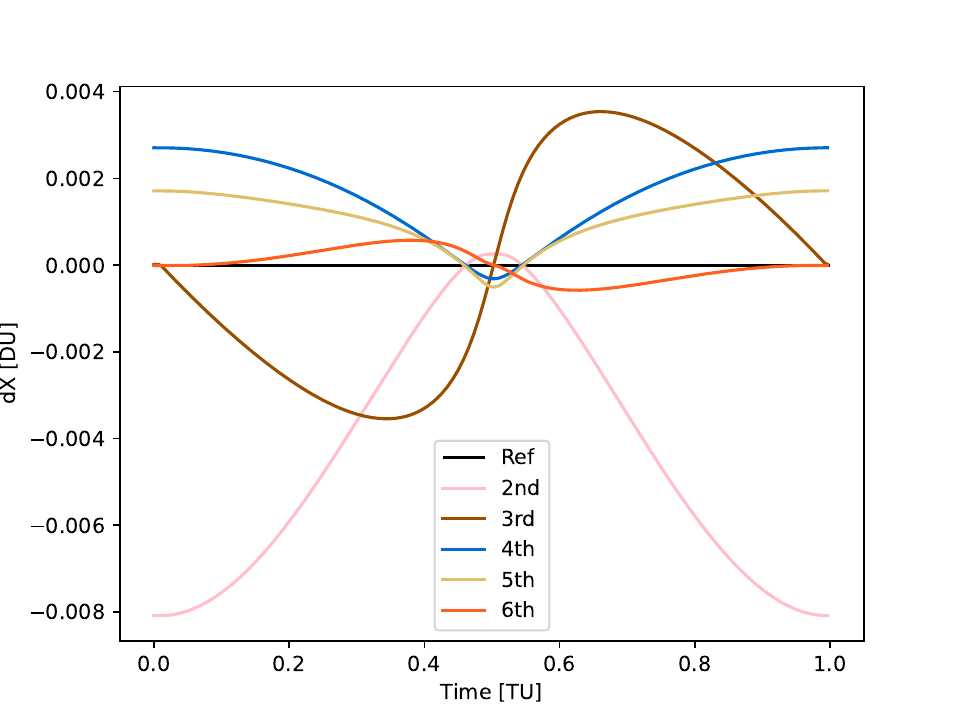}}
  \subfloat{\includegraphics[width=0.5\textwidth]{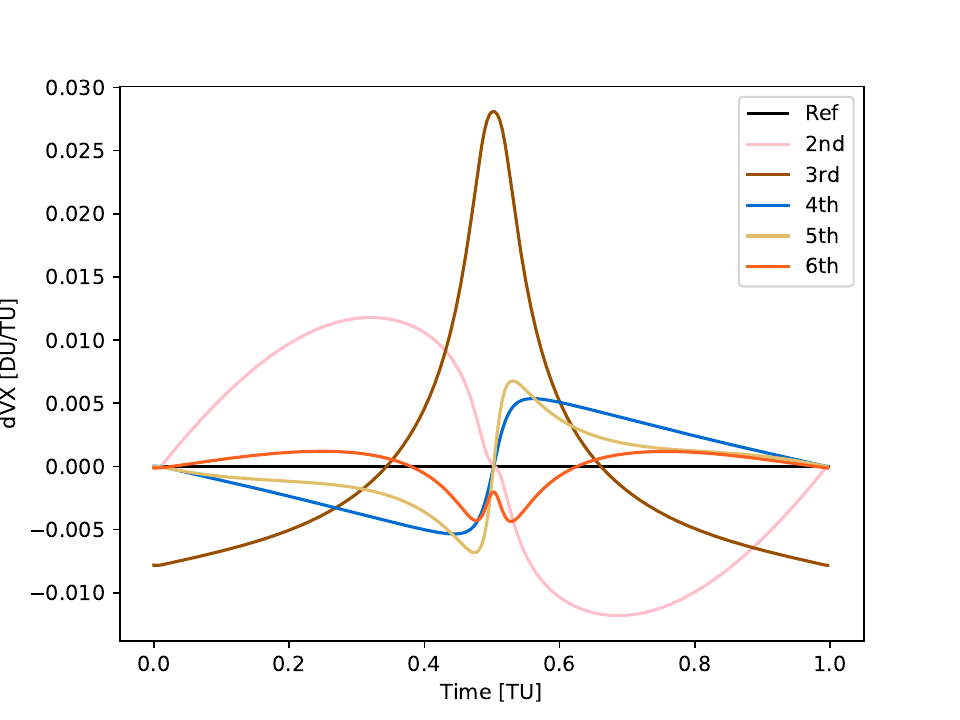}}
  \hfill
  \subfloat{\includegraphics[width=0.5\textwidth]{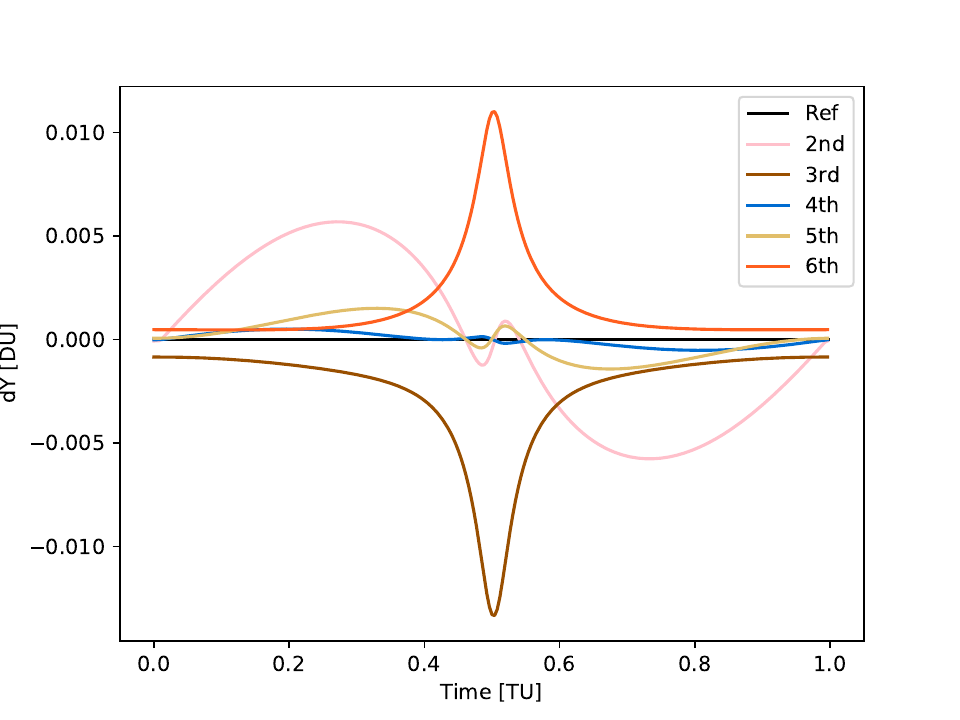}}
  \subfloat{\includegraphics[width=0.5\textwidth]{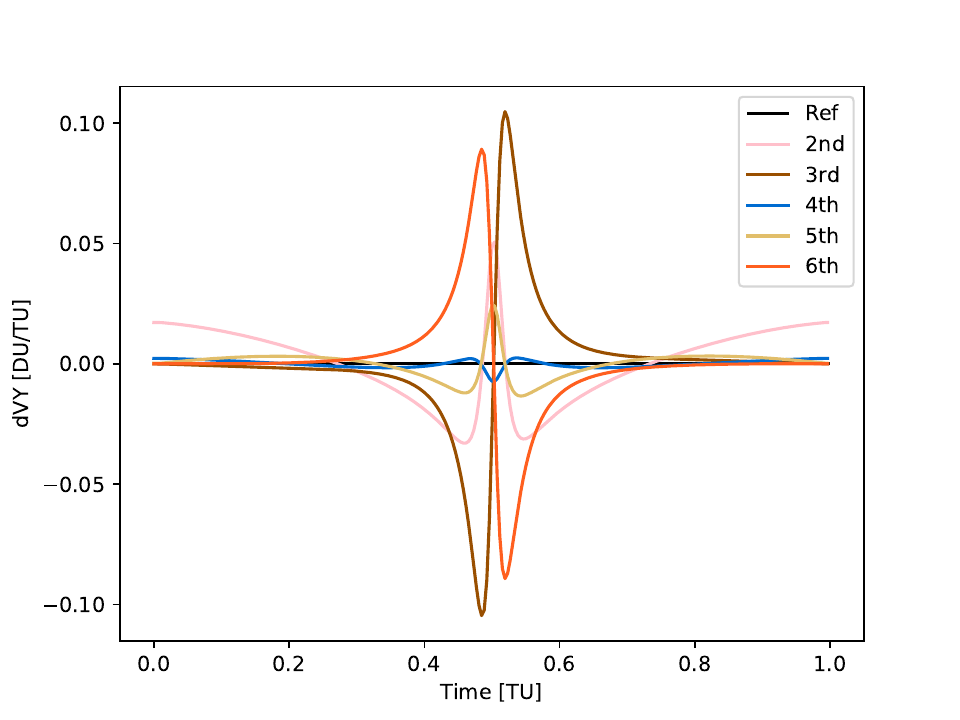}}
  \hfill
  \subfloat{\includegraphics[width=0.5\textwidth]{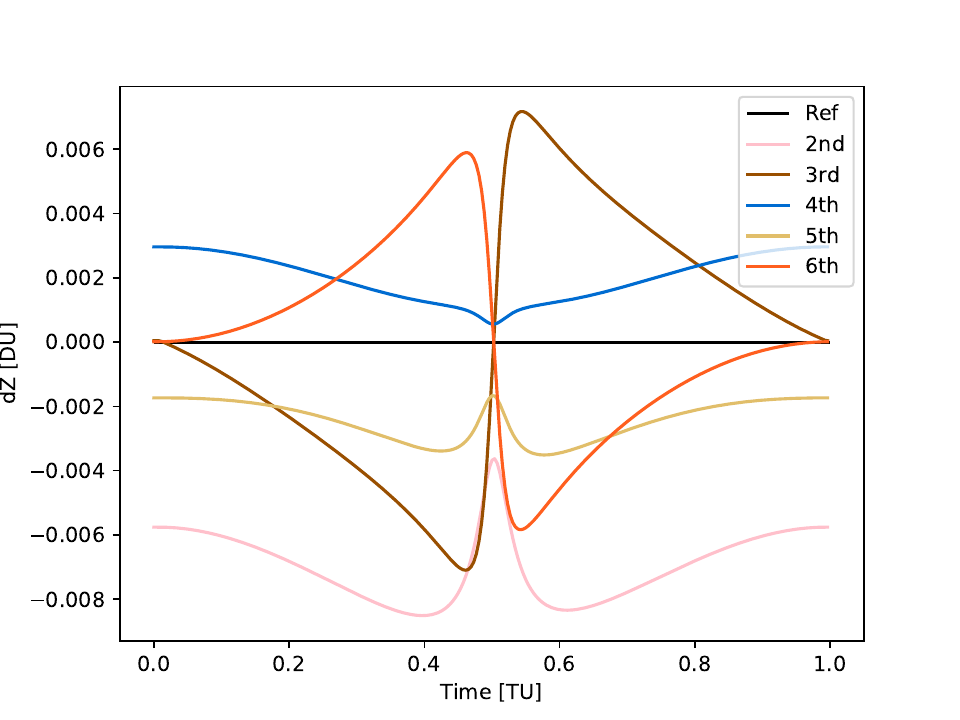}}
  \subfloat{\includegraphics[width=0.5\textwidth]{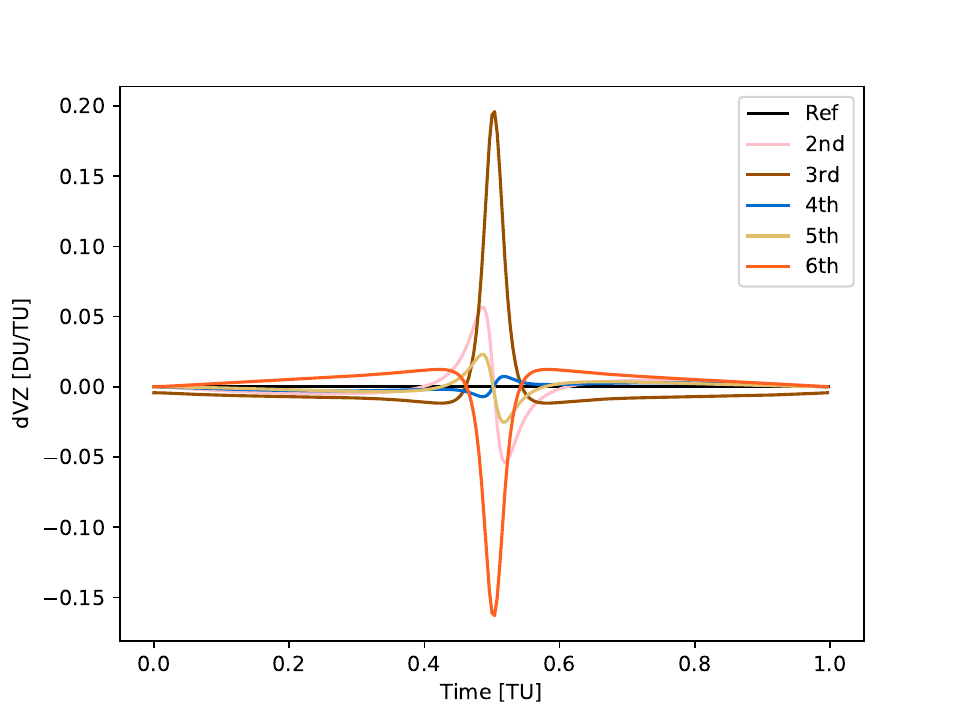}}
  \caption{The difference between the reference trajectory (shown in solid black for all plots) and deviations in the final five eigenvectors. Deviations in the direction of the first eigenvector are not shown, as the extent is theoretically infinite and corresponds to the along-track direction of the orbit.}
    \label{fig:eigen}
\end{figure}

Figure \ref{fig:Trajectories} shows the physical extents in position space of the trajectories. We sample from the boundaries of the six-dimensional reachable set 100,000 times and solve the linearized boundary value problem for each set of initial deviations to obtain initial costates. We then linearly propagate the augmented state deviations and costates and plot all the resulting trajectories in some lower-dimensional space. Due to the constraints of sampling, the reachable bounds in Figure \ref{fig:Trajectories} are not solid and display some gaps and asymmetries. These are products of incomplete sampling and the true reachable set does not display these phenomena. In Figure \ref{fig:Trajectories}, the greatest deviation in the $x$ direction is found at apolune, where the orbit is initialized. This can also be seen in Figure \ref{fig:eigen}, where deviations in the $x$ and $z$ directions are maximized when at apolune. Likewise, the smallest deviations are found at perilune, our half-period. The size of the bounded region is a function of $J^*$, as increasing the cost will expand the reachable set. 

\begin{figure}[ht]
    \centering
    \includegraphics[width=0.5\linewidth]{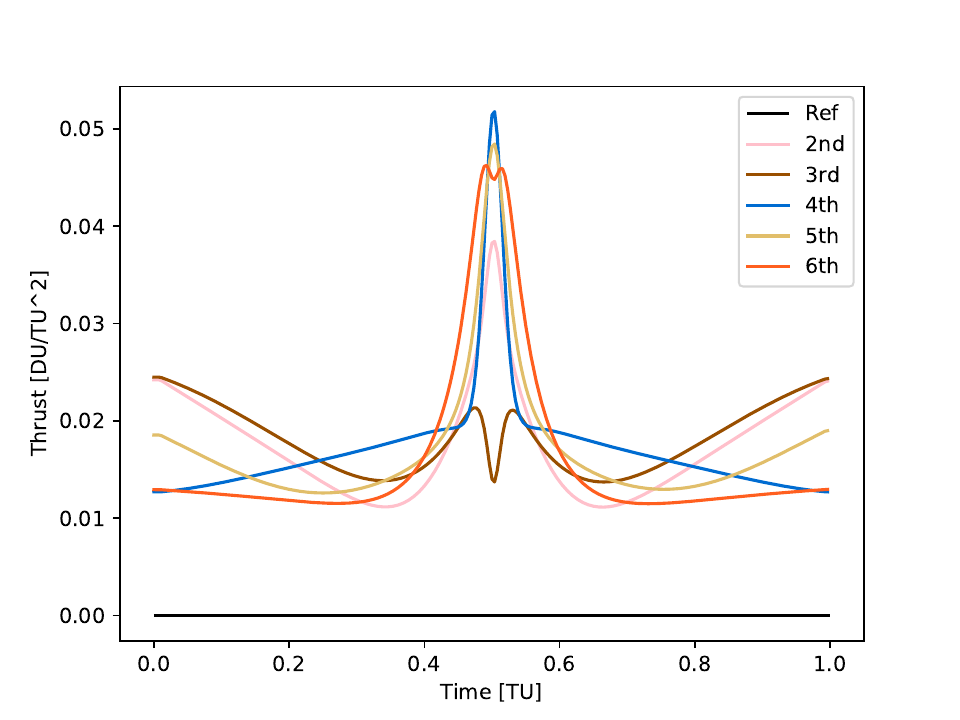}
    \caption{The thrusting magnitude over a single orbit for deviations in the final five eigenvectors.}
    \label{fig:thrust}
\end{figure}

Notice in Figure \ref{fig:eigen} that the orbits begin and conclude at the same point, as expected of periodic trajectories. However, the costates and thus the thrust do not necessarily begin and end at the same values. Thus, the augmented state is not necessarily periodic but the state that strictly includes the position and velocity vectors is periodic.  Figure \ref{fig:eigen} is particularly insightful in how deviations in the initial condition can cascade through the orbit. For example, assume we wish to decrease perilune distance for our spacecraft in order to obtain better images of the Lunar surface. In order to achieve this, the spacecraft should deviate such that $\delta x < 0$, $\delta z < 0$, and $\delta y = 0$ at perilune. From inspection of Figure \ref{fig:eigen}, deviations along the direction of the fifth eigenvector will achieve this, as it fulfills these three qualities. The fifth eigenvector corresponds to the second most-expensive direction to deviate from the reference, so this also indicates that reduction in perilune distance is relatively expensive.


\section{Conclusions} \label{sec:FW}

With this analysis, we have specified bounds on forced periodic trajectories in the CR3BP. These structures are likely to provide unique observation and operational opportunities for satellites. We analyzed these forced periodic trajectories via linear analysis and we determine that for the energy-optimal case, the full augmented state is not periodic but the state that only includes position and velocity is periodic. The linear analysis also allows us to determine the least and most costly directions that deviations to the state may be made. In particular, deviations along the fifth eigenvector listed in Table \ref{tab:eigen}, corresponding to the second-largest eigenvalue, will reduce the perilune distance of a spacecraft for our reference orbit.

\section*{Acknowledgements}
The authors would like to thank Robin Amstrong for helpful discussions on subspace perturbation theory and Jackson Shannon for insightful recommendations. Jackson Kulik would like to thank the DOD National Defense Science and Engineering Graduate Fellowship for support. This material is based upon work supported by the Air Force Office of Scientific Research under award number FA9550-23-1-0665.



\bibliographystyle{AAS_publication}   
\bibliography{references}   

\begin{thebibliography}{10}

\bibitem{Cox2019}
A.~Cox, K.~C. Howell, and D.~C. Folta, ``Dynamical structures in a low-thrust, multi-body model with applications to trajectory design,''  {\em Celestial Mechanics and Dynamical Astronomy}, 2019, 10.1007/s10569-019-9891-7.

\bibitem{Cox2020}
A.~Cox, K.~C. Howell, and D.~C. Folta, ``Trajectory Design Leveraging Low-Thrust, Multi-Body Equilibria and their Manifolds,''  {\em Journal of Astronautical Sciences}, 2020, 10.1007/s40295-020-00211-6.

\bibitem{Morimoto2006}
M.~Morimoto, H.~Yamakawa, and K.~Uesugi, ``Periodic Orbits with Low-Thrust Propulsion in the Restricted Three-body Problem,''  {\em Journal of Guidance, Control, and Dynamics}, Vol.~29, No.~5, 2006, pp.~1131--1139, 10.2514/1.19079.

\bibitem{Tsuruta2024}
A.~Tsuruta, M.~Bando, S.~Hokamoto, and D.~J. Scheeres, ``New Equilibria and Dynamic Structures Under Continuous Optimal Feedback Control,''  {\em Journal of Guidance, Control, and Dynamics}, Vol.~0, No.~0, 2024, pp.~1--12, 10.2514/1.G008270.

\bibitem{DeLeo2022}
L.~DeLeo and M.~Pontani, ``Low-Thrust Orbit Dynamics and Periodic Trajectories in the Earth–Moon System,''  {\em Aerotecnica Missili and Spazio}, 2022, 10.1007/s42496-022-00122-9.

\bibitem{Sandel2024}
C.~Sandel and R.~Sood, ``Natural and Forced Spacecraft Loitering in a Near Rectilinear Halo Orbit,''  {\em Journal of Astronautical Sciences}, Vol.~71, 2024, 10.1007/s40295-024-00446-7.

\bibitem{Kulik2024}
J.~Kulik, M.~Zweig, and D.~Savransky, ``Comparing Relative Reahcable Sets About Nearly Circular Orbits,''  {\em AAS/AIAA Astrodynamics Specialist Conference}, 2024.

\bibitem{bryson2018applied}
A.~E. Bryson, {\em Applied optimal control: optimization, estimation and control}.
\newblock Routledge, 2018.

\bibitem{lee2018reachable}
S.~Lee and I.~Hwang, ``Reachable set computation for spacecraft relative motion with energy-limited low-thrust,''  {\em Aerospace Science and Technology}, Vol.~77, 2018, pp.~180--188.

\bibitem{sun2020analysis}
H.~Sun and J.~Li, ``Analysis on reachable set for spacecraft relative motion under low-thrust,''  {\em Automatica}, Vol.~115, 2020, p.~108864.

\bibitem{Kulik2023}
J.~Kulik, W.~Clark, and D.~Savransky, ``State Transition Tensors for Continuous-Thrust Control of Three-Body Relative Motion,''  {\em Journal of Guidance, Control, and Dynamics}, Vol.~46, No.~8, 2023, pp.~1610--1619, 10.2514/1.G007311.

\end{thebibliography}


\end{document}